\documentclass[iicol,pdflatex,sn-mathphys-num]{sn-jnl}

\usepackage{graphicx}%
\usepackage{multirow}%
\usepackage{amsmath,amssymb,amsfonts}%
\usepackage{amsthm}%
\usepackage{mathrsfs}%
\usepackage[title]{appendix}%
\usepackage{xcolor}%
\usepackage{textcomp}%
\usepackage{manyfoot}%
\usepackage{booktabs}%
\usepackage{algorithm}%
\usepackage{algorithmicx}%
\usepackage{algpseudocode}%
\usepackage{listings}%

\theoremstyle{thmstyleone}%
\theoremstyle{thmstyletwo}%

\theoremstyle{thmstylethree}%

\begin{document}

\title[Stability and Robustness Analysis of Regularized Reconstruction Methods for Low-Dose Computed Tomography in Parallel-Beam Geometry]{Stability and Robustness Analysis of Regularized Reconstruction Methods for Low-Dose Computed Tomography in Parallel-Beam Geometry}


\author[1]{\fnm{Mohamed} \sur{Berrada}}

\affil[1]{\orgdiv{Department of Epidemiology, Public Health and Social Sciences}, \orgname{Faculty of Medicine and Pharmacy of Tangier, Abdelmalek Essaadi University}, 
\orgaddress{
\city{Tangier}, 
\country{Morocco}}}


\abstract{Low-dose computed tomography (LDCT) plays a crucial role in reducing patient radiation exposure, but it increases the ill-posedness of the reconstruction problem due to noise contamination and limited projection data. Regularized reconstruction methods such as Tikhonov and Total Variation (TV) are widely used to improve image quality and stability; however, their performance depends strongly on noise characteristics, sampling conditions, and parameter selection. This study presents a systematic stability and robustness analysis of three classical reconstruction approaches: Filtered Back Projection (FBP), Tikhonov regularization, and TV minimization within a two-dimensional parallel-beam CT framework. A unified simulation pipeline based on the Radon transform is developed and evaluated using both the analytical modified Shepp-Logan phantom and a realistic clinical thorax image. Reconstruction behavior is investigated under multiple degradation scenarios, including Gaussian, Poisson, and mixed noise models, as well as baseline (180 projections) and sparse-view (60 projections) acquisition geometries. To ensure a fair comparison, regularization parameters are optimized for each method and scenario through an exhaustive SSIM-based grid-search procedure. Reconstruction quality is assessed using RMSE, PSNR, and SSIM, while robustness is quantified through an empirical Stability Factor S measuring perturbation amplification from measurement space to image space. The results show that FBP is highly sensitive to noise and angular undersampling. Tikhonov regularization improves structural fidelity compared with FBP but remains more sensitive to perturbation amplification than TV. Conversely, TV provides the best overall compromise between noise suppression, edge preservation, reconstruction accuracy, and numerical stability across both synthetic and anatomical images. These findings highlight the stability-resolution trade-off in LDCT reconstruction and demonstrate that the proposed Stability Factor S provides complementary information to conventional image-quality metrics for assessing reconstruction robustness.
}




\maketitle

\section{Introduction}\label{sec1}

Computed Tomography (CT) is one of the most widely used medical imaging modalities, providing high-resolution cross-sectional visualization of internal anatomical structures for diagnosis, treatment planning, and image-guided interventions \cite{kak2001principles, seeram2015computed}. Despite its considerable clinical value, CT examinations rely on ionizing radiation, motivating continuous efforts to reduce patient exposure while maintaining diagnostically acceptable image quality \cite{brenner2007computed}. This challenge has led to the development of low-dose CT (LDCT) protocols designed to mitigate radiation-induced risks without compromising the clinical utility of reconstructed images  \cite{mccollough2006ct}.

In practice, radiation dose reduction is primarily achieved through two complementary strategies. The first consists of decreasing the X-ray tube current, which reduces the number of detected photons and consequently increases measurement noise \cite{mccollough2009strategies}. The second involves reducing the number of projection views through angular undersampling, thereby lowering radiation exposure and acquisition time at the expense of incomplete measurement information \cite{sidky2008image}. Although effective from a dose-reduction perspective, both approaches substantially degrade reconstruction quality by increasing noise levels, generating streak artifacts, and exacerbating the ill-posed nature of the inverse reconstruction problem \cite{herman2009fundamentals}.

From a mathematical perspective, the objective of CT reconstruction is to recover the cross-sectional attenuation map $f(x,y)$ of the scanned object from its projection measurements collected by a detector array. These measurements are commonly represented in the form of a sinogram and are related to the unknown image through the forward projection operator $A$, yielding the observation model:
\begin{equation}
y = Af + \varepsilon,
\end{equation}
where $y$ denotes the measured projection data and $\varepsilon$ represents measurement perturbations. In low-dose acquisition settings, these perturbations arise from multiple degradation mechanisms. Photon-counting fluctuations due to reduced X-ray flux are naturally described by Poisson statistics, whereas detector electronics introduce additional Gaussian noise components \cite{whiting2006properties}. When combined with angular undersampling, these effects significantly increase reconstruction uncertainty and sensitivity to measurement errors.

The conventional Filtered Backprojection (FBP) algorithm remains the most widely used analytical reconstruction technique owing to its computational efficiency and simplicity \cite{kak2001principles, szczykutowicz2022review}. 
However, FBP is known to be highly sensitive to noise contamination and sparse-view acquisition, often resulting in amplified artifacts and degraded image quality under low-dose conditions. To overcome these limitations, iterative regularized reconstruction methods have been extensively investigated \cite{beister2012iterative}. In a general framework, the reconstruction problem can be formulated as the minimization of the objective function:
\begin{equation}
J(f) = D(Af,y) + \alpha R(f),
\end{equation}
where $D(Af,y)$ is a data-fidelity term, $R(f)$ is a regularization functional encoding prior information, and $\alpha$ is a regularization parameter controlling the trade-off between data fidelity and regularization.

Among the most commonly used regularization approaches, Tikhonov regularization promotes globally smooth solutions by penalizing large fluctuations in the reconstructed image via an $L_2$-norm penalty \cite{tikhonov1977solutions}. In contrast, Total Variation (TV) regularization exploits gradient sparsity via an $L_1$-norm formulation to preserve sharp structural boundaries while suppressing noise and streak artifacts \cite{rudin1992nonlinear, sidky2008image}. Although both methods have demonstrated significant improvements over analytical reconstruction techniques, their effectiveness strongly depends on the acquisition conditions, the noise characteristics, and the selection of the regularization parameter. Furthermore, different regularization strategies may exhibit markedly different trade-offs between reconstruction fidelity, edge preservation, and numerical stability.

Numerous studies have investigated regularized reconstruction methods for low-dose CT and have primarily focused on improving reconstruction accuracy and image quality \cite{sidky2011constrained,loli2022comparison,pan2024compressive}. However, comparatively fewer works have systematically analyzed the robustness of these methods under multiple coupled degradation mechanisms, including photon starvation, electronic noise, and sparse-view acquisition. In particular, the stability of reconstructed images with respect to perturbations in the measurement data remains a critical issue, as small variations in the sinogram may produce large deviations in the reconstructed attenuation map when the inverse problem becomes severely ill-conditioned \cite{arridge2019solving}. Although recent deep learning approaches have shown significant advancements in low-dose CT reconstruction, many contemporary reconstruction frameworks build upon concepts originally developed in classical regularization and model-based reconstruction methods \cite{chen2024deep, leuschner2021lodopab, szczykutowicz2022review}. Consequently, a comprehensive understanding of the stability, robustness, and performance limits of traditional reconstruction methods remains essential for interpreting, auditing, and benchmarking modern data-driven frameworks.

The objective of this work is therefore to provide a comprehensive and reproducible comparative analysis of three representative reconstruction approaches: Filtered Back Projection (FBP), Tikhonov regularization, and Total Variation (TV) regularization within a parallel-beam CT framework. The study investigates their behavior under a broad range of degradation scenarios, including Gaussian noise, Poisson noise, mixed Poisson-Gaussian noise conditions, and angular undersampling configurations involving both full-view (180 projections) and sparse-view (60 projections) acquisitions. To ensure a fair comparison, an exhaustive grid-search strategy is used to determine near-optimal regularization parameters for each scenario, with SSIM serving as the selection criterion. Final performance assessment is conducted using standard image quality metrics, including RMSE, PSNR, and SSIM, together with an empirical Stability Factor S, which quantifies the amplification of measurement perturbations through the reconstruction process. Beyond reconstruction accuracy, particular emphasis is placed on analyzing the balance between structural preservation and numerical stability of each strategy under both an analytical Shepp-Logan phantom and a clinical thorax CT image derived from the LIDC-IDRI database \cite{armato2015data}.

The remainder of this paper is organized as follows. Section 2 reviews the analytical and iterative parallel-beam projection frameworks, providing a  mathematical description of FBP and regularized iterative formulations, and briefly outlining the computational implementation. Section 3 details the experimental setup and evaluation protocol, outlining the phantom configurations, the multi-scenario noise generation, and the exhaustive grid-search framework deployed to optimize hyperparameters. This section also introduces the quantitative structural metrics alongside the relative error amplification factor used to evaluate numerical stability. Finally, Section 4 presents the results and discussion, staging the validation on the analytical Shepp-Logan phantom, the quantitative evaluation on a realistic clinical thorax anatomy obtained from the LIDC-IDRI database, and a comprehensive qualitative visual texture analysis.

\section{Numerical Models and Regularized Iterative Solvers}\label{sec2}
\subsection{Continuous and Discrete Forward Model}\label{sub1sec2}
The physical modeling of X-ray attenuation in parallel-beam Computed Tomography (CT) relies on continuous line integrals through an object $f(x,y)$ \cite{kak2001principles}. Under this geometry, these integrals are parameterized by the projection angle $\theta$ and the perpendicular distance $s$ from the origin. The continuous forward mapping is defined by the Radon transform, $\mathcal{R}$:
\begin{equation}
p(\theta, s) = \mathcal{R}\{f\}(\theta, s)
\end{equation}
where the integral is taken along the line $x \cos\theta + y \sin\theta = s$.

To evaluate the stability and robustness of tomographic reconstruction under low-dose constraints, this continuous formulation is discretized into a finite-dimensional inverse problem. The image $f(x,y)$ is sampled on an $N$-pixel Cartesian grid and reshaped into a vector $\mathbf{f} \in \mathbb{R}^N$, while the sinogram is sampled into a measurement vector $\mathbf{y} \in \mathbb{R}^M$. This leads to the standard linear system:
\begin{equation}
\mathbf{y} = \mathbf{A}\mathbf{f} + \mathbf{e}
\end{equation}
where $\mathbf{A} \in \mathbb{R}^{M \times N}$ represents the discrete projection operator and $\mathbf{e} \in \mathbb{R}^M$ models measurement perturbations. The resulting inverse problem is typically ill-posed, particularly under low-dose conditions, where noise amplification and incomplete angular sampling significantly degrade reconstruction stability \cite{sidky2008image}.

While analytical reconstruction methods remain computationally efficient, their sensitivity to perturbations $\mathbf{e}$ motivates the use of iterative regularized approaches that incorporate prior information into the reconstruction process.
 
\subsection{Analytical Reconstruction: Filtered Backprojection (FBP)}\label{sub2sec2}
The Filtered Backprojection (FBP) algorithm remains the standard analytical reconstruction method in CT due to its computational efficiency and deterministic nature \cite{kak2001principles}. Grounded in the continuous Fourier Slice Theorem, FBP solves the inverse problem by inversion of the Radon transform. Mathematically, the reconstruction of the attenuation profile $f(x,y)$ from the continuous sinogram data $p(\theta, s)$ is formulated as:
\begin{equation}
f(x, y) = \int_{0}^{\pi} \left( p(\theta, s) * q(s) \right) \Big|_{s = x\cos\theta + y\sin\theta} d\theta
\end{equation}
where $*$ denotes the one-dimensional spatial convolution operator, and $q(s)$ represents the spatial representation of the Ram-Lak (ramp) filter. In the frequency domain, this filter scales linearly with absolute frequency, $|\omega|$, acting as a high-pass filter designed to compensate for the $1/r$ blurring intrinsic to simple backprojection.

FBP provides an exact inversion only under ideal conditions of complete angular sampling and noise-free projection data \cite{kak2001principles}. In our discrete algebraic workspace, the continuous backprojection integral is evaluated using the discrete adjoint backprojection operator $\mathbf{A}^T \in \mathbb{R}^{N \times M}$. The discrete FBP pipeline can be compactly expressed as:
\begin{equation}
\mathbf{f}_{\text{FBP}} = \mathbf{A}^T \mathcal{F}^{-1} \left\{ |\omega| \cdot \mathcal{F}\{\mathbf{y}\} \right\}
\end{equation}
where $\mathcal{F}$ and $\mathcal{F}^{-1}$ denote the discrete Fourier transform and its inverse, applied along the detector dimension of the measured sinogram $\mathbf{y}$.

Although computationally efficient, FBP exhibits severe limitations under low-dose constraints. Because the ramp filter $|\omega|$ amplifies high-frequency components, it significantly magnifies the noise energy present in the data. Under clinical low-dose protocols-characterized by photon-starved Poisson noise, detector Gaussian electronic noise, or missing angular views ($\mathbf{e}$)-the high-frequency spectrum of the sinogram is dominated by these perturbations rather than structural information. Consequently, FBP propagates and amplifies this noise into the image space via the adjoint operator $\mathbf{A}^T$, resulting in severe streaking artifacts, substantial image graininess, and a marked loss of low-contrast anatomical details within the parenchymal tissues. These limitations are particularly pronounced in low-dose CT, where noise amplification and angular undersampling severely degrade reconstruction quality \cite{beister2012iterative}.

\subsection{Iterative Regularized Frameworks}\label{sub3sec2}
To overcome the noise amplification and artifacts inherent to analytical inversion under low-dose constraints, iterative regularized approaches reformulate the reconstruction as a mathematical optimization problem \cite{hansen2010discrete,arridge2019solving}. Instead of directly inverting the forward operator, these methods seek an optimal image estimate $\mathbf{f}^*$ by minimizing a generalized objective cost function $J(\mathbf{f})$:
\begin{equation}
\mathbf{f}^* = \arg\min_{\mathbf{f} \geq 0} \left\{ D(\mathbf{Af}, \mathbf{y}) + \alpha R(\mathbf{f}) \right\}
\end{equation}
where $D(\mathbf{Af}, \mathbf{y}) = \frac{1}{2}\|\mathbf{Af} - \mathbf{y}\|_2^2$ represents the $L_2$ squared data-fidelity term enforcing consistency with the measured sinogram $\mathbf{y}$, $R(\mathbf{f})$ is the regularization functional encoding structural a priori knowledge, and $\alpha > 0$ is the regularization parameter. Non-negativity constraints ($\mathbf{f} \geq 0$) are explicitly enforced to reflect the physical reality of X-ray attenuation coefficients.

To systematically evaluate these regularized paradigms against the analytical FBP baseline, this study implements two distinct regularized frameworks within the Core Imaging Library (CIL):
\begin{itemize}
\item \textbf{Tikhonov Regularization ($L_2$-$L_2$ Augmented Framework)}:\\
Tikhonov regularization introduces a smoothness prior by penalizing spatial image variations \cite{tikhonov1977solutions, hansen2010discrete}. The objective function is traditionally formulated as:
\begin{equation}
J_{\text{Tikh}}(\mathbf{f}) = \frac{1}{2}\|\mathbf{Af} - \mathbf{y}\|_2^2 + \frac{\alpha}{2}\|\nabla \mathbf{f}\|_2^2,
\end{equation}
where $\nabla$ represents the discrete spatial gradient operator. 

This formulation can be equivalently rewritten as an augmented least-squares problem:
\begin{equation}
\mathbf{f}^*_{\text{Tikh}} = \arg\min_{\mathbf{f}} \frac{1}{2} \|\mathbf{Kf} - \mathbf{y}_{\text{aug}}\|_2^2,
\end{equation}
where
\begin{equation}
\mathbf{K} = \begin{bmatrix} \mathbf{A} \\ \sqrt{\alpha}\nabla \end{bmatrix} \quad \text{and} \quad \mathbf{y}_{\text{aug}} = [\mathbf{y}, \mathbf{0}]^T,
\end{equation}
with $\mathbf{0}$ denoting a zero-valued vector defined on the range space of the gradient operator.

The resulting quadratic least-squares problem is solved using the Conjugate Gradient Least Squares (CGLS) algorithm \cite{hansen2010discrete} applied directly to the augmented system.

Because the quadratic penalty increases proportionally to the squared gradient magnitude, it promotes globally smooth solutions and suppresses high-frequency variations throughout the image space. While highly effective at mitigating high-frequency Poisson-Gaussian noise, it penalizes genuine sharp transitions just as severely as noise fluctuations, thereby blurring critical diagnostic features such as organ boundaries, vascular structures, and contrast-enhanced interfaces.
\item \textbf{Total Variation Regularization ($L_2$-$L_1$ Primal-Dual Framework)}:\\
Total Variation (TV) regularization addresses the edge-blurring limitations of quadratic regularization by promoting sparsity in the image gradient \cite{rudin1992nonlinear,sidky2008image}. Assuming a piecewise-constant underlying anatomy, the corresponding objective function is defined as

\begin{equation}
J_{\text{TV}}(\mathbf{f})
=
\frac{1}{2}\|\mathbf{Af}-\mathbf{y}\|_2^2
+
\alpha \|\nabla \mathbf{f}\|_{2,1},
\end{equation}

where $\|\nabla \mathbf{f}\|_{2,1}$ denotes the isotropic Total Variation semi-norm, computed as the sum of gradient magnitudes over all image pixels.

Because the TV penalty is non-differentiable, the optimization problem is reformulated into a primal-dual saddle-point framework. Introducing the block operator

\begin{equation}
\mathbf{K}
=
\begin{bmatrix}
\mathbf{A}\\
\nabla
\end{bmatrix},
\end{equation}

the resulting optimization problem is solved using the Primal-Dual Hybrid Gradient (PDHG) algorithm \cite{chambolle2011first}, which alternates between primal and dual updates through proximal mappings associated with the data-fidelity and regularization terms.

Unlike Tikhonov regularization, which penalizes the squared gradient magnitude, TV applies a linear penalty to the gradient norm. Consequently, large sparse gradients corresponding to genuine anatomical boundaries are preserved, while small noise-induced variations are strongly attenuated. This property makes TV particularly effective for low-dose and sparse-view CT reconstruction, although its tendency to promote piecewise-constant solutions may lead to the well-known staircasing artifact.
\end{itemize}

\subsection{Computational Implementation}\label{sub4sec2}

All reconstruction algorithms were implemented within the Core Imaging Library (CIL) framework \cite{j2021core}, which provides a unified environment for tomographic reconstruction and optimization. Forward and backprojection operations were computed using the ASTRA Toolbox \cite{van2015astra}, enabling efficient CPU/GPU-accelerated evaluation of the projection operator ($A$) and its adjoint ($A^T$).

The Tikhonov reconstruction problem, formulated as a quadratic least-squares optimization, was solved using the Conjugate Gradient Least Squares (CGLS) algorithm with a fixed number of 100 iterations. The Total Variation (TV) reconstruction problem was solved using the Primal-Dual Hybrid Gradient (PDHG) algorithm \cite{chambolle2011first}, using 500 iterations to ensure convergence of the non-smooth optimization problem.

\section{Experimental Setup and Evaluation Metrics}
The experimental protocol adopted in this study is organized into four main components. First, we describe the digital phantoms and acquisition geometries used to generate the reference projection data (sinograms). Second, we define the low-dose simulation framework, including the Poisson, Gaussian, and mixed noise models considered throughout the benchmark. Third, we present the hyperparameter optimization strategy used to determine the optimal regularization parameters for each reconstruction method and experimental scenario. Finally, we introduce the quantitative image quality metrics and the empirical stability factor employed to assess reconstruction accuracy, structural preservation, and numerical robustness.

\subsection{Phantoms and Geometry Configuration}
\begin{figure}[htbp]
  \centering
  \includegraphics[width=0.5\textwidth]{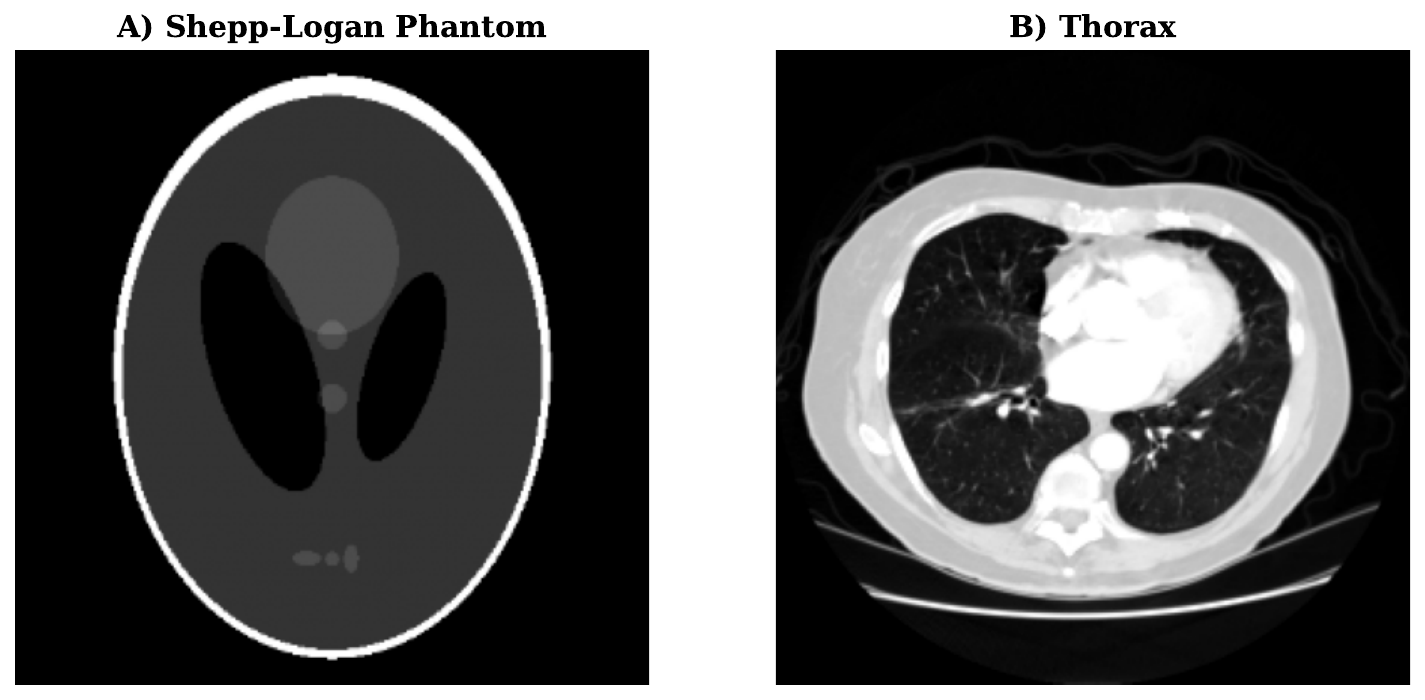}
\caption{Reference images utilized for the computational validation. (A) Analytical modified Shepp--Logan digital phantom (grid size $256 \times 256$ pixels). (B) Clinical thorax CT slice extracted from the LIDC-IDRI database of The Cancer Imaging Archive (TCIA) \cite{armato2015data}, displayed using the default DICOM windowing metadata.} 
\label{fig:FigurePhantomPlusThorax}
\end{figure}

The benchmarking pipeline is executed using two distinct image profiles. The first profile is the analytical, multi-contrast modified Shepp-Logan digital phantom \cite{shepp1974fourier,toft1996radon}, which serves as the standard baseline for testing piecewise-constant approximations and geometric edge preservation. The second profile consists of a real, high-resolution clinical thorax cross-section \cite{armato2015data}, introducing realistic anatomical complexities, continuous parenchymal tissue gradients, and subtle low-contrast variations that challenge the regularized optimization loops (Figure~\ref{fig:FigurePhantomPlusThorax}). Both images are discretized on a uniform Cartesian grid composed of $256\times 256$ pixels.

This $256 \times 256$ pixel matrix establishes the continuous-to-discrete baseline image geometry, where each pixel represents a uniform local attenuation coefficient $\mathbf{f}_{\text{ref}} \in \mathbb{R}^N$ ($N = 65,536$).

The acquisition geometry simulates an ideal parallel-beam CT scanner
 configured via the GPU-accelerated ASTRA-toolbox backend \cite{van2015astra}. The projection space is parameterized by two primary experimental constraints:
\begin{itemize}
\item \textbf{Detector Discretization:}
A linear detector composed of $N_{\text{det}}=256$ bins was employed, matching the image resolution along each spatial dimension. The detector geometry was configured to fully cover the object support, ensuring that no truncation artifacts were introduced during forward projection.
\item \textbf{Angular Sampling:}
Projection angles are uniformly distributed over the interval $[0,\pi)$, corresponding to the complete angular range required for parallel-beam CT acquisition.
\end{itemize}

The pipeline tests two distinct angular sampling densities for each target image to evaluate the robustness of the reconstruction algorithms against structural undersampling. We implement a reference dense-view acquisition consisting of $N_{\theta}=180$ equidistant angular views ($\Delta\theta=1^\circ$) to approximate the best achievable reconstruction quality under the considered acquisition geometry, and a sparse-view acquisition protocol restricted to $N_{\theta}=60$ views ($\Delta\theta=3^\circ$). This undersampled configuration substantially reduces the sinogram measurement dimension from $M = 256 \times 180 = 46\,080$ to $M = 256 \times 60 = 15\,360$, intentionally inducing severe streak artifacts under analytical reconstruction and thereby challenging the edge-preserving capabilities of the regularized iterative solvers.

\subsection{Noise Generation and Low-Dose Scenarios}

We design three distinct experimental noise configurations simulating clinically relevant photon statistics and electronic degradation mechanisms in low-dose CT. Rather than applying uniform perturbations, the benchmark evaluates independent and coupled noise channels across a set of predefined degradation levels, enabling a systematic assessment of reconstruction stability under controlled stochastic conditions.

\begin{itemize}
\item \textbf{Pure Poisson Regime (Photon Starvation):} 
This regime models quantum photon counting statistics under low-dose conditions. For a clean attenuation sinogram $\mathbf{y}^{\mathrm{clean}} = \mathbf{Af}$, the detected photon counts follow:
\begin{equation}
\mathcal{P}_m \sim \mathrm{Poisson}\left(I_0 \exp(-y_m^{\mathrm{clean}})\right),
\end{equation}
where $I_0$ denotes the incident photon flux, varied as $I_0 \in \{10^5, 5\times10^4, 10^4\}$. Lower flux levels simulate severe photon starvation typical of ultra-low-dose CT protocols \cite{whiting2002signal,beister2012iterative,thibault2007three}.

\item \textbf{Pure Gaussian Regime (Electronic Noise):} 
This regime models thermal and readout noise in detector electronics. Additive Gaussian perturbations are applied directly in the sinogram domain:
\begin{equation}
\mathcal{G}_m \sim \mathcal{N}(0,\sigma_G^2),
\quad \text{with} \quad \sigma_G = \eta \cdot \sigma_{\mathrm{true}},
\end{equation}
where $\sigma_{\mathrm{true}}$ is the standard deviation of the noise-free sinogram and $\eta \in \{1\%, 3\%, 5\%\}$ controls the noise intensity.

\item \textbf{Coupled Mixed Regime:} 
This regime represents a realistic acquisition scenario where quantum and electronic noise coexist. The photon flux is fixed to a low-dose setting ($I_0 = 5 \times 10^4$), combined with a Gaussian noise level of $\eta = 3\%$, allowing controlled analysis of coupled degradation effects.
\end{itemize}

For regimes involving quantum noise (Pure Poisson and Mixed), the detected signal is converted into attenuation measurements via a logarithmic transformation. To avoid numerical instabilities due to near-zero photon counts, the detected signal is lower-bounded prior to log compression:
\begin{equation}
y_m = -\ln\left(\frac{\max\left(1,\, \mathcal{P}_m + \mathcal{G}_m\right)}{I_0}\right),
\end{equation}
where $\mathcal{G}_m = 0$ in the Pure Poisson case and $\mathcal{G}_m \sim \mathcal{N}(0,\sigma_G^2)$ in the Mixed regime. This non-linear transformation induces a heteroscedastic and signal-dependent noise structure in the sinogram domain, significantly increasing the ill-posedness of the inverse problem.

The experiments consider seven noise scenarios ordered by increasing severity: Pure Gaussian (1\%), Pure Poisson ($10^5$), Pure Poisson ($5 \times 10^4$), Pure Gaussian (3\%), Mixed (3\%, $5 \times 10^4$), Pure Poisson ($10^4$), and Pure Gaussian (5\%).

\subsection{Hyperparameter Optimization Framework}

A key challenge in regularized iterative tomographic reconstruction is the selection of the regularization parameter $\alpha$, which controls the trade-off between data fidelity and the regularization penalty \cite{hansen2010discrete}. In this benchmark, we adopt an exhaustive grid search to determine optimal parameters for each reconstruction model and experimental configuration.

The optimal regularization parameters $(\alpha_{\text{TV}}^{\text{opt}}, \alpha_{\text{Tikh}}^{\text{opt}})$ are obtained by maximizing the Structural Similarity Index (SSIM) with respect to the ground-truth image:
\begin{equation}
\alpha_{\text{opt}} = \arg\max_{\alpha \in \Lambda} \text{SSIM}(\mathbf{f}_\alpha^*, \mathbf{f}_{\text{ref}})
\end{equation}
where $\mathbf{f}_\alpha^*$ denotes the reconstructed image for a given $\alpha$, and the search space is defined as:
\begin{equation}
\Lambda = \{10^k \, | \, k \in [-5,1], \Delta k = 0.25\}
\end{equation}

This logarithmic grid ensures a uniform exploration across several orders of magnitude, covering both over-regularized and under-regularized regimes.

The optimization is performed independently for each phantom type (Shepp-Logan and clinical thorax) due to their distinct structural characteristics. The Shepp-Logan phantom exhibits piecewise-constant regions, favoring higher regularization levels, whereas the clinical thorax contains smooth anatomical variations requiring more conservative regularization to avoid excessive smoothing and staircasing artifacts. The complete list of optimal hyperparameter values is provided in Appendix \ref{ann:hyperparameter}, Table~\ref{tab:optimized_hyperparameters}.

\subsection{Quantitative Evaluation Metrics}

To evaluate reconstruction accuracy and structural fidelity, we compute three standard image quality metrics against the ground-truth image $\mathbf{f}_{\text{ref}}$.

\begin{enumerate}
\item Root Mean Squared Error (RMSE)\\
The RMSE measures pixel-wise intensity deviations between the reconstruction $\mathbf{f}^*$ and the reference $\mathbf{f}_{\text{ref}}$:
\begin{equation}
\text{RMSE}(\mathbf{f}^*, \mathbf{f}_{\text{ref}}) = \sqrt{\frac{1}{N} \sum_{j=1}^{N} (f^*_j - f_{\text{ref},j})^2}
\end{equation}
Lower values indicate higher numerical accuracy.
\item Peak Signal-to-Noise Ratio (PSNR)\\
The PSNR provides a logarithmic measure of reconstruction quality:
\begin{equation}
\text{PSNR}(\mathbf{f}^*, \mathbf{f}_{\text{ref}}) =
10 \log_{10} \left( \frac{L^2}{\text{MSE}} \right)
\end{equation}
where $L$ denotes the dynamic range of the reference image, and $\text{MSE} = \text{RMSE}^2$ represents the Mean Squared Error.

Higher values indicate better noise suppression and artifact reduction.
\item Structural Similarity Index Measure (SSIM)\\
The SSIM evaluates perceptual similarity between $x = \mathbf{f}^*$ and $y = \mathbf{f}_{\text{ref}}$ \cite{wang2004image}:
\begin{equation}
\text{SSIM}(x,y) =
\frac{(2\mu_x\mu_y + C_1)(2\sigma_{xy} + C_2)}
{(\mu_x^2 + \mu_y^2 + C_1)(\sigma_x^2 + \sigma_y^2 + C_2)}
\end{equation}
where $\mu_x$, $\mu_y$ are local means, $\sigma_x^2$, $\sigma_y^2$ local variances, and $\sigma_{xy}$ the cross-covariance. The constants $C_1$ and $C_2$ are included to avoid numerical instability. The SSIM index is bounded between $-1$ and $1$, with $1$ indicating perfect structural agreement.
\end{enumerate}

\paragraph{Stability Analysis via Error Amplification.}
While the above metrics quantify reconstruction fidelity in the image domain, they do not explicitly characterize the sensitivity of the inverse problem to perturbations in the measurement space. This aspect is closely related to classical stability analysis in inverse problems \cite{engl1996regularization}.

To this end, we define the Relative Error Amplification Factor:
\begin{equation}
S =
\frac{\|\mathbf{f}_{\delta} - \mathbf{f}_{\text{ref}}\|_2 / \|\mathbf{f}_{\text{ref}}\|_2}
{\|\mathbf{y}_{\delta} - \mathbf{y}_{\text{ref}}\|_2 / \|\mathbf{y}_{\text{ref}}\|_2}
\end{equation}
where $\mathbf{f}_{\text{ref}}$ and $\mathbf{y}_{\text{ref}}$ denote the reference image and sinogram, and $\mathbf{f}_{\delta}$ and $\mathbf{y}_{\delta}$ their noise-perturbed counterparts under identical acquisition settings.

The quantity $S$ measures how perturbations in the measurement domain are amplified through the reconstruction pipeline. Values close to $S \approx 1$ indicate stable behavior, whereas larger values indicate increased sensitivity to noise and under-sampling. Importantly, $S$ should be interpreted as an empirical stability indicator rather than a formal condition number, since the reconstruction process includes discretization, regularization, and iterative nonlinear solvers.

\section{Results and Discussion}\label{sec:results}

This section presents the quantitative and qualitative evaluation of the three reconstruction frameworks considered in this study: Filtered Backprojection (FBP), Tikhonov regularization, and Total Variation (TV) regularization. The methods are assessed on both the modified Shepp--Logan phantom and a clinical thorax image under seven degradation scenarios combining Gaussian, Poisson, and mixed Poisson--Gaussian noise models.

To investigate the impact of angular undersampling, all experiments are performed using both a full-view acquisition geometry (180 projections) and a sparse-view configuration (60 projections). Reconstruction quality is evaluated using RMSE, PSNR, and SSIM, while numerical robustness is assessed through the Relative Error Amplification Factor ($S$), which quantifies the propagation of measurement perturbations into the reconstructed image domain.

The results are first analyzed on the analytical Shepp--Logan phantom to characterize the intrinsic behavior of the reconstruction algorithms under controlled conditions, and subsequently on the clinical thorax image to assess their performance in a more realistic anatomical setting.

\subsection{Results on the Shepp--Logan Phantom}
The modified Shepp--Logan phantom provides a controlled benchmark for assessing the behavior of reconstruction algorithms in the presence of noise and angular undersampling. Owing to its piecewise-constant structure and well-defined boundaries, this phantom is particularly suitable for evaluating edge preservation, artifact suppression, and numerical stability under extreme acquisition conditions. 

\begin{figure}[htbp]
  \centering
  \includegraphics[width=0.5\textwidth]{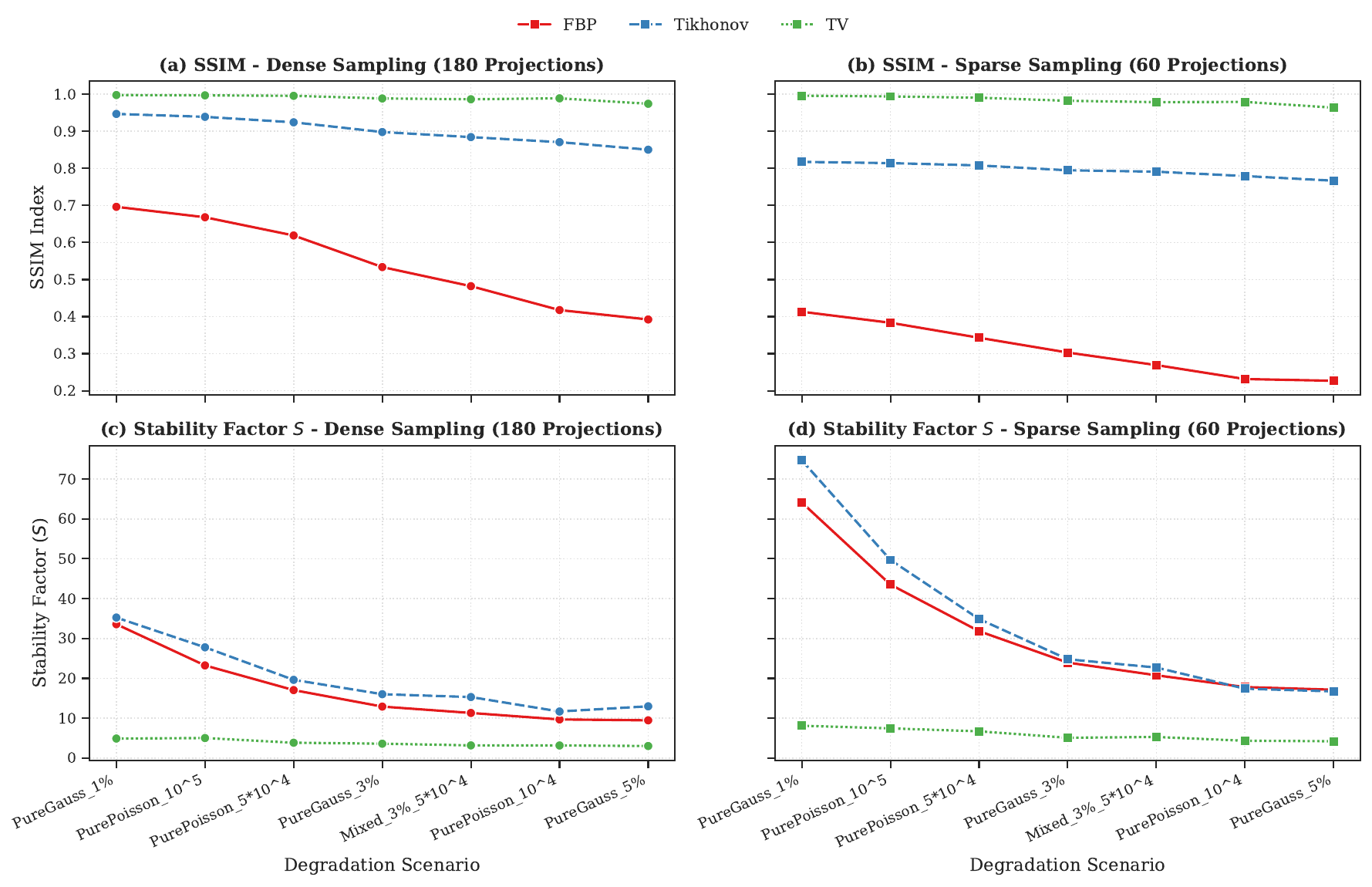}
  \caption{Scenario-by-scenario quantitative evaluation of the Structural Similarity Index (SSIM) and the empirical Stability Factor ($S$) on the modified Shepp--Logan phantom. The curves track performance across the seven distinct degradation levels for: (a) SSIM under baseline conditions (180 projections), (b) SSIM under sparse-view constraints (60 projections), (c) Stability Factor $S$ under baseline conditions (180 projections), and (d) Stability Factor $S$ under sparse-view constraints (60 projections).}
 \label{fig:figure_metrics_comparison_phantom}
\end{figure}

\begin{figure*}[htbp]
  \centering
  \includegraphics[width=0.9\textwidth]{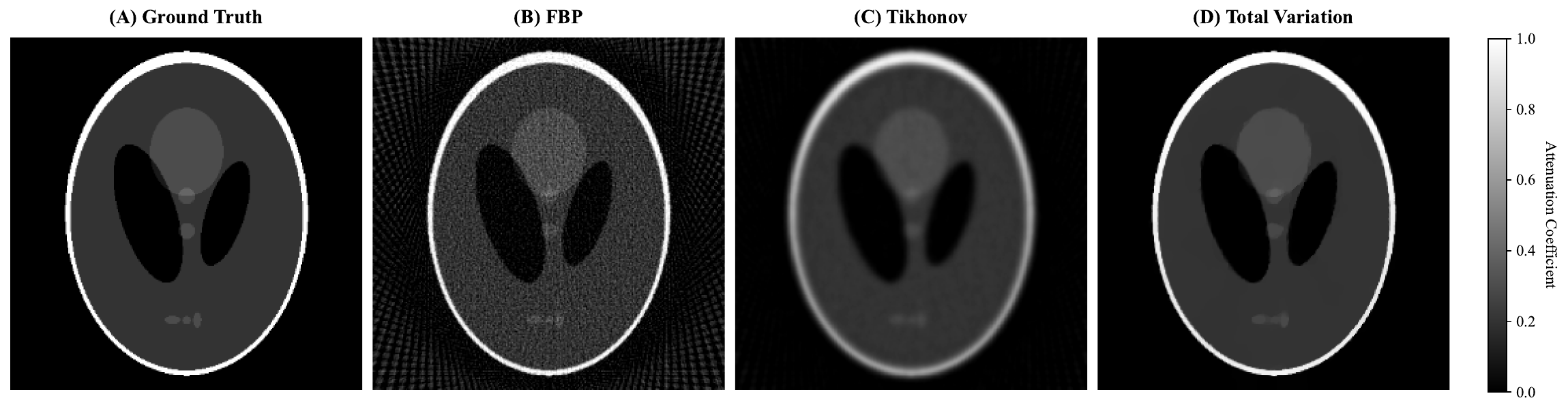}
  \caption{Representative visual reconstruction results for the modified Shepp--Logan phantom under a critical sparse-view and mixed noise acquisition scenario ($60$ projections, coupled Mixed noise regime with $\eta = 3\%$ Gaussian noise and $I_0 = 5 \times 10^4$ photon flux). From left to right: Ground Truth, Filtered Backprojection (FBP), Tikhonov regularization, and Total Variation (TV) regularization.}
  \label{fig:visual_comparison_phantom_60views_Mixed_G3_P5x10p4}
\end{figure*}

Overall, the quantitative and qualitative results reveal a clear hierarchy between the reconstruction methods. Table~\ref{tab:mean_std_phantom} summarizes the average quantitative performance metrics (RMSE, PSNR, SSIM) and the empirical stability factor ($S$) of FBP, Tikhonov, and TV reconstructions across the seven degradation scenarios. The detailed scenario-by-scenario numerical results are reported in Appendix~\ref{ann:hyperparameter}, Table~\ref{tab:results_phantom}, for exhaustiveness. Complementary graphical analyses displaying the scenario-dependent evolution of the Structural Similarity Index (SSIM) and the Stability Factor ($S$) are illustrated in Figure~\ref{fig:figure_metrics_comparison_phantom}, while representative reconstructed visual fields are displayed in Figure~\ref{fig:visual_comparison_phantom_60views_Mixed_G3_P5x10p4}.

As compiled in Table~\ref{tab:mean_std_phantom}, TV consistently provided the best quantitative performance across both acquisition geometries. Under the reference 180-view configuration, TV achieved the lowest RMSE, highest PSNR and SSIM values, together with the lowest stability factor $S$, confirming its superior reconstruction performance even under relatively favorable sampling conditions. This advantage became increasingly evident under the sparse-view configuration (60 views), where TV achieved an average SSIM of 0.984 compared with 0.796 for Tikhonov and 0.310 for FBP, while reducing the average stability factor from approximately 31--34 for FBP and Tikhonov to 5.890. 

This behavior highlights the robustness of the total variation prior, whose sparse-gradient assumption is particularly well adapted to the piecewise-constant structure of the Shepp--Logan phantom. In contrast, FBP exhibited a strong degradation in both image fidelity and stability, reaching an average stability factor of 31.306 in the 60-view configuration, indicating substantial sensitivity to measurement perturbations. Tikhonov regularization improved structural similarity compared with FBP, but its quadratic smoothing prior resulted in limited robustness under sparse-view conditions, with an average stability factor of 34.414, remaining considerably higher than that of TV (5.890).

\begin{table*}[t]
\caption{Overall quantitative reconstruction performance on the modified Shepp--Logan phantom. Metrics represent the mean $\pm$ standard deviation computed across the seven distinct degradation scenarios (comprising pure Poisson, pure Gaussian, and coupled mixed noise regimes) for both the baseline (180 views) and sparse-view (60 views) acquisition geometries. Bold values indicate the top-performing method for each configuration.}
\label{tab:mean_std_phantom}
\centering
\small
\begin{tabular}{llcccc}
\toprule
\textbf{Geometry} & \textbf{Method} & \textbf{RMSE} & \textbf{PSNR (dB)} & \textbf{SSIM} & $\bf S$  \\
\midrule
180 views & FBP & $0.0432\pm 0.006$ & $27.359\pm 1.195$ & $0.544\pm 0.121$ & $16.753\pm 8.873$ \\
          & Tikhonov & $0.0533\pm 0.011$ & $25.636\pm 1.865$ & $0.902\pm 0.036$ & $19.809\pm 8.651$  \\
          & TV & $\bf0.0114\pm 0.004$ & $\bf 39.450\pm 3.552$ & $\bf 0.990\pm 0.008$ & $\bf 3.820\pm 0.830$  \\
\midrule
60 views  & FBP & $0.0800\pm 0.010$ & $21.990\pm 1.063$ & $0.310\pm 0.073$ & $31.306\pm 17.189$  \\
          & Tikhonov & $0.0850\pm 0.005$ & $21.429\pm 0.541$ & $0.796\pm 0.019$ & $34.414\pm 21.168$  \\
          & TV & $\bf0.0171\pm 0.005$ & $\bf 35.783\pm 3.110$ & $\bf 0.984\pm 0.011$ & $\bf 5.890\pm 1.538$  \\
 \bottomrule
\end{tabular}
\end{table*}

Figure~\ref{fig:figure_metrics_comparison_phantom} provides a detailed scenario-by-scenario analysis of structural fidelity and numerical stability. While all evaluated methods experienced an expected performance degradation as noise severity increased, TV regularization effectively decoupled reconstruction accuracy from noise magnitude, maintaining SSIM values close to unity and exhibiting significantly lower error amplification factors than both FBP and Tikhonov. The SSIM evolution of Tikhonov regularization remained substantially closer to TV than to FBP, indicating improved structural preservation through quadratic smoothing. However, its stability factor remained closer to that of FBP, highlighting the limitation of purely quadratic regularization in controlling perturbation amplification in severely ill-conditioned reconstruction problems. This observation emphasizes that high structural similarity does not necessarily translate into enhanced numerical robustness, as image fidelity and perturbation sensitivity characterize complementary aspects of reconstruction performance. The stability advantage of TV became particularly evident under the sparse-view geometry (60 views), where FBP suffered from pronounced streak artifacts and noise amplification, while Tikhonov reduced these artifacts at the cost of increased smoothing and loss of fine structural details.

Visual inspections of the reconstructed slices confirm these quantitative findings. As illustrated in Figure~\ref{fig:visual_comparison_phantom_60views_Mixed_G3_P5x10p4}, FBP reconstructions suffer from highly damaging streak artifacts and non-stationary noise contamination, which completely mask internal structures under sparse-view conditions. Tikhonov regularization isotropic smoothing effectively eliminates these streak patterns but introduces severe edge blurring, penalizing the visibility of small high-contrast features. In contrast, TV reconstruction successfully suppresses noise and streaking while maintaining sharp, crisp object boundaries. The reconstructed fields remain visually close to the ground truth, validating the structural enforcement of the gradient sparsity model even under coupled mixed degradation regimes.

Overall, the Shepp-Logan benchmark confirms that Total Variation regularization provides the best compromise between noise suppression, edge preservation, and reconstruction stability. The advantage of TV becomes increasingly pronounced as both noise intensity and angular undersampling increase, highlighting the suitability of sparse-gradient priors for piecewise-constant structures and low-dose CT acquisition protocols.

This superiority is partly expected because the modified Shepp-Logan phantom closely satisfies the piecewise-constant assumption underlying TV regularization. Consequently, the phantom represents a favorable scenario for sparse-gradient models and provides an ideal benchmark for evaluating edge-preserving reconstruction behavior.

\subsection{Results on the Clinical Thorax Image}

While the Shepp--Logan phantom provides a useful benchmark for evaluating reconstruction performance on piecewise-constant structures, it does not fully reflect the complexity of clinical CT imaging. The thorax image introduces realistic anatomical characteristics, including continuous tissue density variations, low-contrast pulmonary structures, vascular networks, and heterogeneous attenuation patterns. Consequently, this clinical test case constitutes a substantially more challenging testbed for assessing the practical robustness of reconstruction algorithms under low-dose acquisition conditions.

\begin{figure}[htbp]
  \centering
  \includegraphics[width=0.5\textwidth]{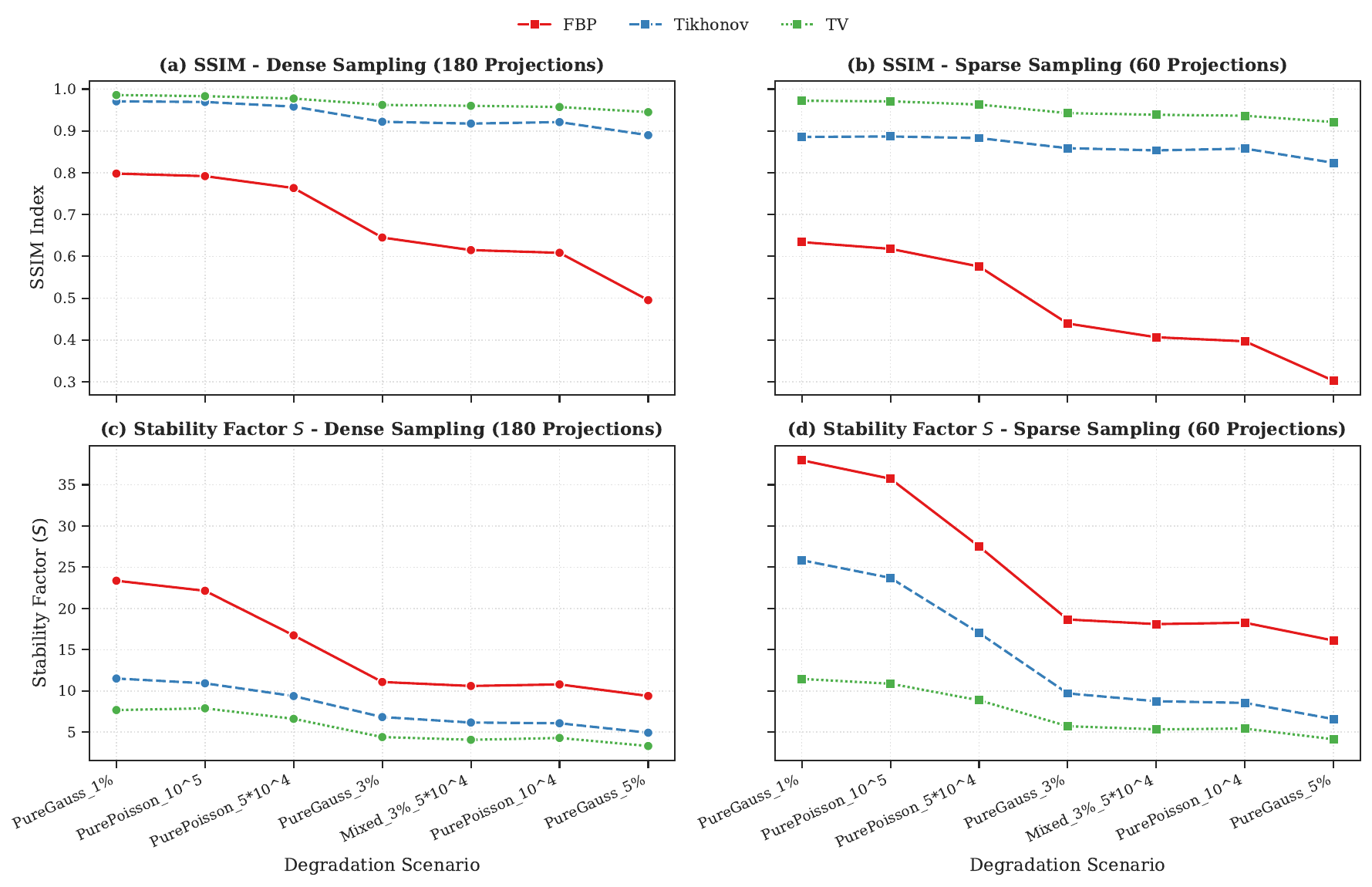}
\caption{Scenario-by-scenario quantitative evaluation of the Structural Similarity Index (SSIM) and the empirical Stability Factor ($S$) on the clinical thorax CT image. The curves track performance across the seven distinct degradation levels for: (a) SSIM under baseline conditions (180 projections), (b) SSIM under sparse-view constraints (60 projections), (c) Stability Factor $S$ under baseline conditions (180 projections), and (d) Stability Factor $S$ under sparse-view constraints (60 projections).}  \label{fig:figure_metrics_comparison_thorax}
\end{figure}

\begin{figure*}[htbp]
  \centering
  \includegraphics[width=0.9\textwidth]{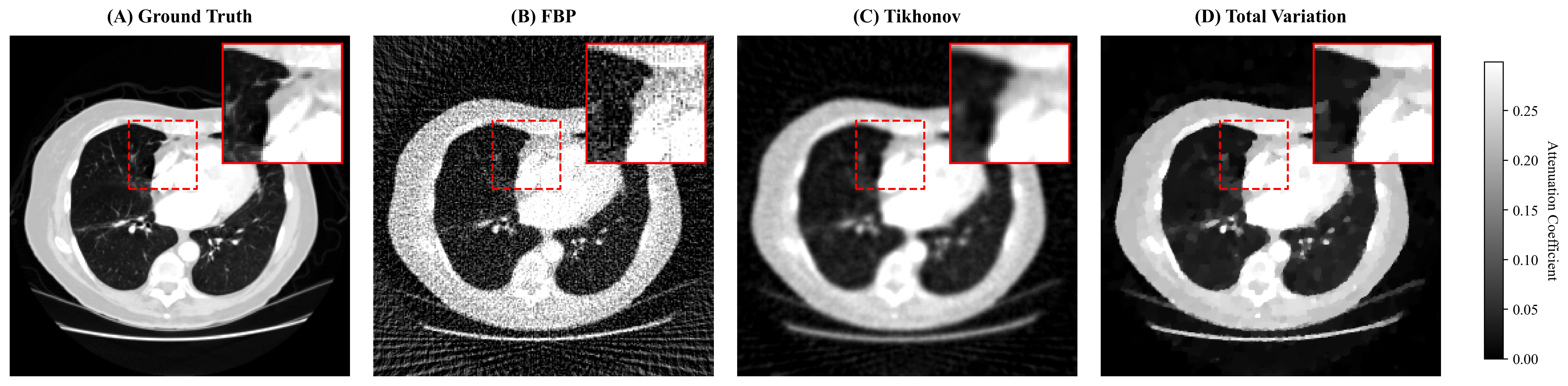}
\caption{Representative visual reconstruction results for the clinical thorax CT image under a critical sparse-view and mixed noise acquisition scenario ($60$ projections, coupled Mixed noise regime with $\eta = 3\%$ Gaussian noise and $I_0 = 5 \times 10^4$ photon flux). From left to right: Ground Truth, Filtered Backprojection (FBP), Tikhonov regularization, and Total Variation (TV) regularization. The inset in the upper-right corner of each panel provides a magnified view of a selected region of interest, highlighting the structural preservation of fine anatomical textures and edge-sharpening capabilities of each framework under severe degradation.}  \label{fig:visual_comparison_zoom_thorax_60views_Mixed_G3_P5x10p4}
\end{figure*}

\begin{figure*}[htbp]
  \centering
  \includegraphics[width=0.9\textwidth]{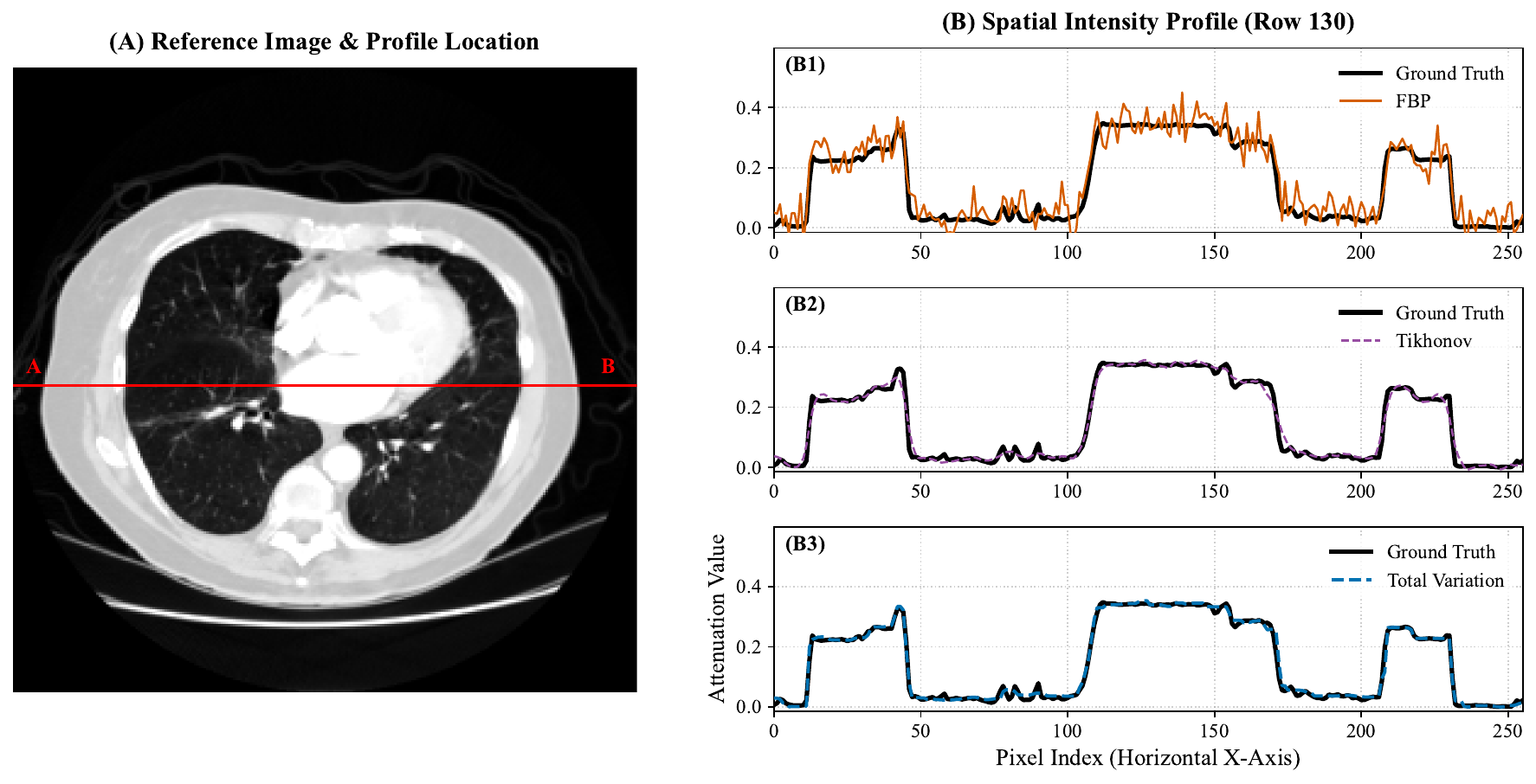}
\caption{Spatial intensity profile analysis on the clinical thorax CT image under severe sparse-view and mixed noise degradation ($60$ projections). (A) Reference image showing the spatial location of the extracted horizontal cross-section (Row 130). (B) Corresponding 1D spatial intensity profiles comparing each framework against the Ground Truth (GT): (B1) Filtered Backprojection (FBP) vs. GT, (B2) Tikhonov regularization vs. GT, and (B3) Total Variation (TV) regularization vs. GT.}  \label{fig:Line_Profile_Row130_Mixed_G3_P5x10p4}
\end{figure*}

Table~\ref{tab:mean_std_thorax} summarizes the average reconstruction performance obtained across the seven degradation scenarios for both the 180-view and 60-view acquisition geometries. Detailed scenario-specific numerical results are provided in Appendix~\ref{ann:hyperparameter}, Table~\ref{tab:results_thorax}. The evolution of the Structural Similarity Index (SSIM) and the Stability Factor ($S$) across all degradation conditions is illustrated in Figure~\ref{fig:figure_metrics_comparison_thorax}, while representative reconstruction examples with magnified regions of interest are displayed in Figure~\ref{fig:visual_comparison_zoom_thorax_60views_Mixed_G3_P5x10p4}. Additional intensity line profiles are presented in Figure~\ref{fig:Line_Profile_Row130_Mixed_G3_P5x10p4} to facilitate a localized assessment of edge preservation and contrast fidelity.

\begin{table*}[t]
\caption{Overall quantitative reconstruction performance on the clinical thorax CT image. Metrics represent the mean $\pm$ standard deviation computed across the seven distinct degradation scenarios (comprising pure Poisson, pure Gaussian, and coupled mixed noise regimes) for both the baseline (180 views) and sparse-view (60 views) acquisition geometries. Bold values indicate the top-performing method for each configuration.}
\label{tab:mean_std_thorax}
\centering
\small
\begin{tabular}{llcccccccc}
\toprule
\textbf{Geometry} & \textbf{Method} & \textbf{RMSE} & \textbf{PSNR (dB)} & \textbf{SSIM} & $\bf S$  \\
\midrule
180 views & FBP & $0.024 \pm 0.006$ & $31.962 \pm 2.248$ & $0.674 \pm 0.114$ & $14.877 \pm 5.884$ \\
 & Tikhonov & $0.013 \pm 0.004$ & $37.216 \pm 2.656$ & $0.935 \pm 0.031$ & $7.980 \pm 2.600$  \\
 & TV & $\bf 0.009 \pm 0.003$ & $\bf 40.531 \pm 2.524$ & $\bf 0.967 \pm 0.015$ & $\bf 5.475 \pm 1.877$\\
\midrule
60 views & FBP & $0.041 \pm 0.012$ & $27.558 \pm 2.417$ & $0.482 \pm 0.127$ & $24.628 \pm 9.136$  \\
 & Tikhonov & $0.021 \pm 0.002$ & $32.904 \pm 0.747$ & $0.864 \pm 0.023$ & $14.315 \pm 7.894$ \\
 & TV & $\bf 0.012 \pm 0.003$ & $\bf 38.092 \pm 1.997$ & $\bf 0.949 \pm 0.019$ & $\bf 7.410 \pm 2.949$  \\
 \bottomrule
\end{tabular}
\end{table*}

Unlike the Shepp--Logan phantom, where the piecewise-constant assumption naturally favors sparse-gradient regularization, the clinical thorax image contains extended smooth transitions and subtle low-contrast anatomical structures. Nevertheless, TV regularization remained the overall best-performing reconstruction framework across both acquisition geometries. Averaged over all degradation scenarios, TV consistently achieved the highest structural similarity and the lowest perturbation amplification factors. Under the 180-view configuration, TV reached an average SSIM of 0.967 and a stability factor of 5.475, compared with 0.935 and 7.980 for Tikhonov and 0.674 and 14.877 for FBP, respectively. This demonstrates that TV preserves anatomical information while providing improved numerical robustness even when the underlying image does not strictly satisfy the sparse-gradient assumption.

The quantitative differences between the reconstruction methods become particularly evident in the sparse-view geometry. When reducing the number of projections from 180 to 60, FBP exhibited a pronounced deterioration in reconstruction quality, with the SSIM decreasing from 0.674 to 0.482 and the stability factor increasing from 14.877 to 24.628, reflecting enhanced sensitivity to angular undersampling artifacts. Tikhonov regularization successfully suppressed part of these artifacts and improved structural fidelity compared with FBP, achieving an SSIM of 0.864 and a stability factor of 14.315. However, this improvement was obtained through stronger spatial smoothing, which limits the preservation of fine anatomical details. In contrast, TV maintained the best compromise between artifact suppression and structural preservation, achieving an SSIM of 0.949 and the lowest stability factor (7.410) under the challenging sparse-view configuration.

Figure~\ref{fig:figure_metrics_comparison_thorax} confirms these observations. Although all methods exhibit a gradual decrease in SSIM as noise severity increases, the evolution of the stability factor $S$ reflects a different behavior, since it quantifies the relative amplification of measurement perturbations rather than the absolute reconstruction error. TV exhibits the most favorable overall stability profile, maintaining the lowest amplification factors across the investigated scenarios while preserving the highest structural fidelity. Unlike the Shepp--Logan phantom, the gap between TV and Tikhonov is reduced for the clinical thorax anatomy, as the heterogeneous anatomical structures and smooth intensity transitions partially satisfy the assumptions of quadratic regularization. Nevertheless, TV provides the most consistent compromise between artifact suppression, structural preservation, and robustness under degraded acquisition conditions.

The visual comparisons presented in Figure~\ref{fig:visual_comparison_zoom_thorax_60views_Mixed_G3_P5x10p4} further highlight the qualitative differences between reconstruction strategies, with the magnified insets providing a closer inspection of fine anatomical structures and tissue interfaces. FBP reconstructions suffer from strong streaking artifacts that obscure local structures and introduce significant intensity distortions, particularly under the sparse-view configuration. Tikhonov regularization effectively attenuates these artifacts but produces noticeable over-smoothing, reducing the visibility of fine details and sharp transitions. In contrast, TV regularization preserves the delineation of structural boundaries within the zoomed regions while maintaining a substantially cleaner background appearance, demonstrating its ability to balance artifact suppression and structural preservation under challenging low-dose acquisition conditions.

To further investigate local reconstruction fidelity, Figure~\ref{fig:Line_Profile_Row130_Mixed_G3_P5x10p4} displays the 1D spatial intensity profiles extracted along Row 130. The FBP profile (Panel B1) exhibits strong high-frequency oscillations and intensity fluctuations, directly reflecting the structured streak artifacts introduced by sparse-view and noisy acquisition conditions. While Tikhonov regularization (Panel B2) effectively suppresses these fluctuations, it fails to accurately reproduce abrupt intensity transitions, leading to excessive smoothing at structural interfaces and reduced local contrast preservation. In contrast, TV regularization (Panel B3) demonstrates the highest agreement with the Ground Truth profile, accurately preserving steep gradients associated with anatomical boundaries. However, a slight stair-casing effect remains observable in regions characterized by slowly varying intensity patterns, reflecting the intrinsic limitation of gradient-sparsity priors for smoothly varying structures.

Overall, the clinical thorax experiments demonstrate that the benefits of TV regularization extend beyond idealized piecewise-constant phantoms and remain evident under realistic anatomical imaging conditions. Although the performance gap with respect to Tikhonov reconstruction is reduced compared with that observed for the Shepp--Logan phantom, TV consistently provides the most favorable compromise between quantitative accuracy, structural preservation, and numerical stability, particularly under sparse-view and low-dose acquisition conditions.

\subsection{Comparative Analysis}

The joint evaluation of both synthetic and clinical images highlights critical insights into the interplay between regularization priors, anatomical morphology, and numerical stability.

A key finding is that the performance advantage of TV regularization over Tikhonov regularization was considerably more pronounced for the analytical Shepp--Logan phantom than for the clinical thorax image. This behavior is consistent with the mathematical formulations underlying these regularization strategies. The Shepp--Logan phantom is inherently piecewise-constant, meaning that its spatial gradient is highly sparse, which closely aligns with the sparse-gradient prior promoted by the TV functional. Conversely, the realistic clinical thorax image exhibits complex structural textures, fine anatomical details, and continuous intensity transitions that only partially satisfy the piecewise-constant assumption. While TV still achieved the highest SSIM values and the lowest stability factor $S$ for the thorax image, this morphology-dependent variation highlights that the effectiveness of gradient-sparsity regularization strongly depends on the underlying anatomical complexity.

Across both imaging configurations, angular undersampling represented a major degradation factor, particularly when comparisons were performed under identical noise conditions. Reducing the number of projection views from 180 to 60 resulted in a consistent decrease in SSIM and an increase in the stability factor $S$ across the evaluated reconstruction methods. This observation indicates that incomplete angular information constitutes a critical source of instability in the investigated low-dose CT configurations. While statistical noise mainly introduces intensity fluctuations that can be partially controlled through regularization, sparse-view acquisitions generate structural ambiguities and streak artifacts by increasing the number of poorly constrained components of the inverse problem. Consequently, stronger prior information, such as the gradient sparsity enforced by TV, becomes essential to compensate for the missing angular information.

Furthermore, the proposed empirical Stability Factor $S$ revealed vulnerability trends that were not fully captured by conventional image-quality metrics alone. While RMSE and PSNR quantify global reconstruction fidelity, they do not directly characterize the sensitivity of the reconstruction process to measurement perturbations. Consequently, two reconstruction methods exhibiting similar image-quality indicators may still present substantially different stability behaviors, which becomes apparent through the analysis of $S$.

Quantitatively, TV regularization consistently exhibited the most favorable stability behavior, maintaining values of $S$ below 10 even under mixed-noise and sparse-view conditions. In comparison, FBP and Tikhonov reconstructions showed higher sensitivity to perturbations, particularly under severe angular undersampling, with stability factors substantially larger than those obtained with TV. This demonstrates that enforcing an $L_1$ gradient prior not only improves structural boundary preservation but also enhances robustness against stochastic perturbations introduced in the projection data.

In conclusion, this comparative analysis demonstrates that TV regularization provides the most favorable compromise between reconstruction accuracy, structural preservation, and numerical stability across the investigated acquisition conditions. Although its advantages are particularly pronounced for piecewise-constant structures, TV remains highly effective for realistic anatomical imaging scenarios involving sparse-view and low-dose conditions. Finally, the Stability Factor $S$ provides an additional perspective beyond conventional image-quality metrics, enabling a more comprehensive assessment of the practical robustness of tomographic reconstruction methods.

It should be noted that the present study was restricted to two-dimensional parallel-beam geometries and classical variational regularization models. Although these settings provide a controlled framework for benchmarking, clinical CT systems typically involve more complex three-dimensional acquisition geometries and reconstruction pipelines. Extending the proposed stability analysis to three-dimensional acquisitions, higher-order regularization strategies, and emerging deep-learning-based reconstruction methods represents a promising direction for future research.

\section{Conclusion}\label{sec:conclusion}
This study presented a systematic benchmark of three tomographic reconstruction paradigms-Filtered Backprojection (FBP), Tikhonov regularization, and Total Variation (TV) minimization under realistic low-dose acquisition conditions combining sparse-view sampling, Poisson quantum noise, Gaussian electronic noise, and mixed degradation regimes. The evaluation was conducted on both an analytical modified Shepp-Logan phantom and a realistic clinical thorax image to assess reconstruction accuracy, structural preservation, and numerical stability. By optimizing the regularization parameters independently for each acquisition scenario, this benchmarking framework ensured a fair comparison of the reconstruction methods under identical acquisition and optimization conditions.

The experimental results consistently demonstrated that TV regularization provided the best overall compromise across all investigated acquisition settings, achieving the highest structural similarity, the lowest reconstruction errors, and the strongest robustness against both noise and angular undersampling. While this advantage was particularly pronounced for the piecewise-constant Shepp-Logan phantom whose sparse-gradient structure closely matches the assumptions underlying the TV prior-TV remained the most effective framework for the anatomically complex thorax image. 

Beyond conventional image-quality assessment, this work introduced the empirical Stability Factor ($S$) as a complementary metric to explicitly quantify the propagation and amplification of measurement perturbations through the inverse reconstruction process. Under this stability perspective, TV regularization consistently exhibited the lowest amplification factors, indicating a markedly reduced sensitivity to measurement perturbations compared with FBP and Tikhonov regularization. These findings provide additional insight into the resolution-stability trade-off that governs tomographic reconstruction under low-dose conditions.

Although the present study focused on two-dimensional parallel-beam geometries and classical variational models, the proposed benchmarking methodology and stability analysis are sufficiently general to be extended to more complex acquisition settings. Future work will investigate three-dimensional cone-beam CT configurations, advanced sparsity-promoting regularizers, and emerging deep-learning-based reconstruction approaches. Ultimately, the proposed stability framework could provide a valuable and reproducible tool for assessing the numerical robustness of modern hybrid and data-driven reconstruction algorithms under clinically realistic low-dose conditions.

\bmhead{Acknowledgements}
The author would like to acknowledge the Core Imaging Library (CIL) and ASTRA-toolbox developers for providing the open-source computational frameworks utilized in this benchmarking study.

\bmhead{Code Availability}
The complete tomographic reconstruction pipeline developed in this work is publicly available under the project name LODESTAR. The source code is hosted on GitHub (\url{https://github.com/MBerrada-FMPT/LODESTAR}) and is permanently archived on Zenodo \cite{berrada_lodestar_2026}.

\bmhead{Author Contributions}
M.B. Conceptualization, Methodology, Software, Validation, Formal analysis, Investigation, Writing -- Original Draft.

\section*{Declarations}
\textbf{Conflict of interest} The author declare no competing interests.

\begin{appendices}

\section{Optimal hyperparameters and Detailed Quantitative Results}\label{ann:hyperparameter}

\begin{sidewaystable}
\caption{Scenario-specific optimal hyperparameters ($\alpha_{\text{opt}}$) maximizing the Structural Similarity Index (SSIM) evaluated under fully-sampled ($180$ views) and under-sampled ($60$ views) tomographic geometries.}\label{tab:optimized_hyperparameters}
\begin{tabular*}{\textwidth}{@{\extracolsep\fill}@{} lllcccc@{}}
\toprule
& & & \multicolumn{2}{c}{$\alpha_{\text{TV}}$} & \multicolumn{2}{c}{$\alpha_{\text{Tikh}}$} \\
\cmidrule(lr){4-5} \cmidrule(lr){6-7}
Phantom Type & Degradation Channel & Noise Level & 180 views & 60 views & 180 views & 60 views \\
\midrule
\multirow{7}{*}{Shepp-Logan}     & \multirow{3}{*}{Pure Gaussian} & $\sigma = 1\%$ & $3.16\times10^{-3}$ & $1.78\times10^{-3}$ & $0.01$ & $0.0562$ \\
                                 &                                & $\sigma = 3\%$ & $0.01$ & $5.62\times10^{-3}$ & $0.0316$ & $0.0562$ \\
                                 &                                & $\sigma = 5\%$ & $0.0178$ & $0.01$ & $0.1$ & $0.1$ \\
                                 \cmidrule(lr){2-7}
                                 & \multirow{3}{*}{Pure Poisson}  & $I_0 = 10^5$    & $5.62\times10^{-3}$ & $3.16\times10^{-3}$ & $0.0178$ & $0.0562$ \\
                                 &                                & $I_0 = 5\times 10^4$ & $5.62\times10^{-3}$ & $5.62\times10^{-3}$ & $0.0178$ & $0.0562$ \\
                                 &                                & $I_0 = 10^4$    & $0.0178$ & $0.01$ & $0.0562$ & $0.1$ \\
                                 \cmidrule(lr){2-7}
                                 & Mixed                          & $\sigma = 3\%,\; I_0=5\times 10^4$ & $0.01$ & $0.01$ & $0.0562$ & $0.1$ \\
\midrule
\multirow{7}{*}{Clinical Thorax} & \multirow{3}{*}{Pure Gaussian} & $\sigma = 1\%$ & $1.00\times10^{-3}$ & $5.62\times10^{-4}$ & $5.62\times 10^{-3}$ & $0.0178$ \\
                                 &                                & $\sigma = 3\%$ & $5.62\times10^{-3}$ & $3.16\times10^{-3}$ & $0.0316$ & $0.0316$ \\
                                 &                                & $\sigma = 5\%$ & $0.01$ & $5.62\times10^{-3}$ & $0.0562$ & $0.0562$ \\
                                 \cmidrule(lr){2-7}
                                 & \multirow{3}{*}{Pure Poisson}  & $I_0 = 10^5$    & $1.00\times10^{-3}$ & $5.62\times10^{-4}$ & $5.62\times 10^{-3}$ & $0.0178$ \\
                                 &                                & $I_0 = 5\times 10^4$ & $1.78\times10^{-3}$ & $1.00\times10^{-3}$ & $0.01$ & $0.0178$ \\
                                 &                                & $I_0 = 10^4$    & $5.62\times10^{-3}$ & $3.16\times10^{-3}$ & $0.0316$ & $0.0316$ \\
                                 \cmidrule(lr){2-7}
                                 & Mixed                          & $\sigma = 3\%,\; I_0=5\times 10^4$ & $5.62\times10^{-3}$ & $3.16\times10^{-3}$ & $0.0316$ & $0.0316$ \\
\bottomrule
\end{tabular*}
\end{sidewaystable}


\begin{sidewaystable*}
\caption{Exhaustive quantitative performance metrics (RMSE, PSNR, SSIM, and exact Stability Factor $S$) for the analytical Shepp-Logan phantom under Full-View (180 projections) and Sparse-View (60 projections) reconstruction regimes across the seven standardized degradation scenarios.}
\label{tab:results_phantom}
\centering
\small
\begin{tabular}{llcccccccc}
\toprule
 & & \multicolumn{4}{c}{\textbf{Full-View Geometry (180 Projections)}} & \multicolumn{4}{c}{\textbf{Sparse-View Geometry (60 Projections)}} \\
\cmidrule(lr){3-6} \cmidrule(lr){7-10}
\textbf{Scenario} & \textbf{Method} & \textbf{RMSE} & \textbf{PSNR (dB)} & \textbf{SSIM} & \textbf{$S$} & \textbf{RMSE} & \textbf{PSNR (dB)} & \textbf{SSIM} & \textbf{$S$} \\
\midrule
\texttt{PureGauss 1\%} & FBP & 0.0367 & 28.71 & 0.6964 & 33.58 & 0.0692 & 23.20 & 0.4132 & 64.11 \\
 & Tikhonov & 0.0385 & 28.29 & 0.9471 & 35.25 & 0.0806 & 21.87 & 0.8176 & 74.70 \\
 & TV & \textbf{0.0054} & \textbf{45.43} & \textbf{0.9981} & \textbf{4.90} & \textbf{0.0088} & \textbf{41.16} & \textbf{0.9963} & \textbf{8.11} \\
\midrule
\texttt{PurePoisson $10^5$} & FBP & 0.0375 & 28.53 & 0.6682 & 23.24 & 0.0706 & 23.02 & 0.3834 & 43.52 \\
 & Tikhonov & 0.0448 & 26.98 & 0.9393 & 27.78 & 0.0807 & 21.87 & 0.8143 & 49.69 \\
 & TV & \textbf{0.0081} & \textbf{41.83} & \textbf{0.9972} & \textbf{5.03} & \textbf{0.0121} & \textbf{38.35} & \textbf{0.9943} & \textbf{7.45} \\
\midrule
\texttt{PurePoisson $5\cdot10^4$} & FBP & 0.0391 & 28.16 & 0.6190 & 17.05 & 0.0736 & 22.66 & 0.3430 & 31.83 \\
 & Tikhonov & 0.0450 & 26.94 & 0.9244 & 19.62 & 0.0807 & 21.86 & 0.8081 & 34.88 \\
 & TV & \textbf{0.0088} & \textbf{41.07} & \textbf{0.9959} & \textbf{3.86} & \textbf{0.0155} & \textbf{36.18} & \textbf{0.9908} & \textbf{6.71} \\
\midrule
\texttt{PureGauss 3\%} & FBP & 0.0423 & 27.48 & 0.5338 & 12.92 & 0.0781 & 22.15 & 0.3031 & 23.95 \\
 & Tikhonov & 0.0524 & 25.61 & 0.8982 & 16.02 & 0.0808 & 21.85 & 0.7950 & 24.79 \\
 & TV & \textbf{0.0118} & \textbf{38.55} & \textbf{0.9888} & \textbf{3.61} & \textbf{0.0166} & \textbf{35.59} & \textbf{0.9826} & \textbf{5.10} \\
\midrule
\texttt{Mixed 3\% + $5\cdot10^4$} & FBP & 0.0451 & 26.92 & 0.4822 & 11.32 & 0.0829 & 21.63 & 0.2692 & 20.77 \\
 & Tikhonov & 0.0610 & 24.30 & 0.8846 & 15.31 & 0.0907 & 20.85 & 0.7912 & 22.71 \\
 & TV & \textbf{0.0127} & \textbf{37.94} & \textbf{0.9867} & \textbf{3.18} & \textbf{0.0211} & \textbf{33.49} & \textbf{0.9790} & \textbf{5.30} \\
\midrule
\texttt{PurePoisson $10^4$} & FBP & 0.0505 & 25.94 & 0.4177 & 9.69 & 0.0925 & 20.67 & 0.2315 & 17.81 \\
 & Tikhonov & 0.0609 & 24.31 & 0.8709 & 11.69 & 0.0904 & 20.88 & 0.7794 & 17.40 \\
 & TV & \textbf{0.0165} & \textbf{35.67} & \textbf{0.9892} & \textbf{3.16} & \textbf{0.0226} & \textbf{32.94} & \textbf{0.9796} & \textbf{4.34} \\
\midrule
\texttt{PureGauss 5\%} & FBP & 0.0514 & 25.78 & 0.3923 & 9.47 & 0.0933 & 20.60 & 0.2268 & 17.15 \\
 & Tikhonov & 0.0705 & 23.03 & 0.8505 & 12.99 & 0.0909 & 20.83 & 0.7669 & 16.72 \\
 & TV & \textbf{0.0165} & \textbf{35.67} & \textbf{0.9745} & \textbf{3.03} & \textbf{0.0230} & \textbf{32.78} & \textbf{0.9640} & \textbf{4.22} \\
\bottomrule
\end{tabular}
\end{sidewaystable*}


\begin{sidewaystable*}
\caption{Exhaustive quantitative performance metrics (RMSE, PSNR, SSIM, and exact Stability Factor $S$) for the clinical thorax CT image under Full-View (180 projections) and Sparse-View (60 projections) reconstruction regimes across the seven standardized degradation scenarios.}
\label{tab:results_thorax}
\centering
\small
\begin{tabular}{llcccccccc}
\toprule
 & & \multicolumn{4}{c}{\textbf{Full-View Geometry (180 Projections)}} & \multicolumn{4}{c}{\textbf{Sparse-View Geometry (60 Projections)}} \\
\cmidrule(lr){3-6} \cmidrule(lr){7-10}
\textbf{Scenario} & \textbf{Method} & \textbf{RMSE} & \textbf{PSNR (dB)} & \textbf{SSIM} & \textbf{$S$} & \textbf{RMSE} & \textbf{PSNR (dB)} & \textbf{SSIM} & \textbf{$S$} \\
\midrule
PureGauss 1\% & FBP & 0.0178 & 34.46 & 0.7979 & 23.37 &  0.029 & 30.290 & 0.634 & 37.974
 \\
              & Tikhonov & 0.0088 & 40.61 & 0.9706 & 11.51 & 0.020 & 33.629 & 0.886 & 25.854 \\
              & TV & \textbf{0.0059} & \textbf{44.12} & \textbf{0.9857} & \textbf{7.69} & \textbf{0.009} & \textbf{40.701} & \textbf{0.972} & \textbf{11.454} \\
\midrule
PurePoisson $10^5$ & FBP & 0.0181 & 34.32 & 0.7918 & 22.15 & 0.030 & 30.070 & 0.618 & 35.734 \\
                   & Tikhonov & 0.0089 & 40.46 & 0.9691 & 10.93 & 0.020 & 33.632 & 0.887 & 23.714 \\
                   & TV & \textbf{0.0065} & \textbf{43.28} & \textbf{0.9830} & \textbf{7.90} & \textbf{0.009} & \textbf{40.396} & \textbf{0.971} & \textbf{10.884} \\
\midrule
PurePoisson $5\cdot10^4$ & FBP & 0.0195 & 33.68 & 0.7633 & 16.74 & 0.032 & 29.402 & 0.576 & 27.528 \\
                         & Tikhonov & 0.0109 & 38.71 & 0.9580 & 9.38 & 0.020 & 33.572 & 0.883 & 17.030 \\
                         & TV & \textbf{0.0077} & \textbf{41.74} & \textbf{0.9774} & \textbf{6.62} & \textbf{0.010} & \textbf{39.237} & \textbf{0.963} & \textbf{8.872} \\
\midrule
PureGauss 3\% & FBP & 0.0253 & 31.43 & 0.6450 & 11.09 & 0.042 & 26.954 & 0.439 & 18.667 \\
              & Tikhonov & 0.0156 & 35.63 & 0.9218 & 6.84 & 0.022 & 32.635 & 0.859 & 9.705\\
              & TV & \textbf{0.0101} & \textbf{39.44} & \textbf{0.9621} & \textbf{4.41} & \textbf{0.013} & \textbf{37.218} & \textbf{0.942} & \textbf{5.726} \\
\midrule
Mixed 3\% + $5\cdot10^4$ & FBP & 0.0271 & 30.82 & 0.6149 & 10.61 & 0.046 & 26.308 & 0.407 & 18.110 \\
                        & Tikhonov & 0.0158 & 35.53 & 0.9175 & 6.17 & 0.022 & 32.618 & 0.853 & 8.758 \\
                        & TV & \textbf{0.0104} & \textbf{39.12} & \textbf{0.9600} & \textbf{4.08} & \textbf{0.013} & \textbf{36.902} & \textbf{0.939} & \textbf{5.348} \\
\midrule
PurePoisson $10^4$ & FBP & 0.0282 & 30.48 & 0.6084 & 10.80 & 0.047 & 26.004 & 0.397 & $18.274$ \\
                   & Tikhonov & 0.0159 & 35.46 & 0.9211 & 6.09 & 0.022 & 32.599 & 0.858 & 8.552 \\
                   & TV & \textbf{0.0112} & \textbf{38.47} & \textbf{0.9571} & \textbf{4.30} &  \textbf{0.014} & \textbf{36.527} & \textbf{0.936} & \textbf{5.441} \\
\midrule
PureGauss 5\% & FBP & 0.0353 & 28.53 & 0.4953 & 9.39 & 0.060 & 23.878 & 0.303 & 16.111 \\
              & Tikhonov & 0.0186 & 34.11 & 0.8898 & 4.94 &  0.025 & 31.641 & 0.824 & 6.591 \\
              & TV & \textbf{0.0125} & \textbf{37.54} & \textbf{0.9449} & \textbf{3.33} &  \textbf{0.016} & \textbf{35.667} & \textbf{0.921} & \textbf{4.147} \\
\bottomrule
\end{tabular}
\end{sidewaystable*}

\end{appendices}



\end{document}